\documentclass[preprint,11pt]{article}%
\usepackage{graphicx}
\usepackage{amsmath,amsxtra,amssymb,latexsym, amscd,amsthm}
\usepackage[mathscr]{eucal}
\newtheorem{thm}{Theorem}[section]

\newtheorem{lem}{Lemma}[section]

\theoremstyle{definition}

\newtheorem{exam}{Example}[section]

\theoremstyle{remark}
\newtheorem{rem}{Remark}[section]
\numberwithin{equation}{section}

\begin{document}

\title%{Gaussian concentration for the vector of convex functions}
%{A variance-sensitive concentration inequality for the vector of non-Lipschitz functions}
{A variance-sensitive Gaussian concentration inequality}
\author{Nguyen Tien Dung
\thanks{Department of Mathematics, VNU University of Science, Vietnam National University, Hanoi, 334 Nguyen
Trai, Thanh Xuan, Hanoi, 084 Vietnam. Email: dung@hus.edu.vn}\,\,\footnote{Department of Mathematics, FPT University, Hoa Lac High Tech Park, Hanoi, Vietnam.}}%Email: dungnt@fpt.edu.vn

\date{\today}          % Ngay

\maketitle
\begin{abstract} In this note, we obtain a Gaussian concentration inequality for a class of non-Lipschitz functions. In the one-dimensional case, our results supplement those established by Paouris and Valettas  in \cite{Paouris2018}.
\end{abstract}
%\section{Introduction}       % Muc dau tien
%{\large Densities, Tail probabilities and }
\section{Introduction}
In the whole paper, let $Z=(Z_1,\cdots,Z_n)$ be an $n$-dimensional standard Gaussian vector and $Z'=(Z'_1,\cdots,Z'_n)$ be an independent copy of $Z.$ It is well known that if $f:\mathbb R^n \to \mathbb R$ is a Lipschitz function with Lipschitz constant $L,$ then we have the following concentration inequalities (see e.g. \cite{Bogachev1998,Pisier1986}), for all $t>0,$
\begin{equation}\label{aa1}
P\left(f(Z)-Ef(Z)  \leq-t\right)\leq \exp\left(-\frac{t^2}{2L^2}\right)
\end{equation}
and 
\begin{equation}\label{aba1}
P\left(f(Z)-M  \leq -t\right)\leq \frac{1}{2}\exp\left(-\frac{t^2}{2L^2}\right),
\end{equation}
where $M$ is a median for $f (Z).$ In a recent paper \cite{Paouris2018}, Paouris and Valettas pointed out that Gaussian concentration phenomenon also holds for the convex functions. More specifically, they got the following.

\noindent{\bf Theorem A.} {\it Let $f:\mathbb R^n \to \mathbb R$ be a convex function with $E|f(Z)|<\infty$ and let $M$ be a median for $f (Z).$ Then, it holds that
\begin{equation}\label{a2d3}
P\left(f(Z)-M  \leq-tE(f(Z)-M)_+\right)\leq \Phi\left(-\frac{\sqrt{2\pi}}{32}t\right)\,\,\,\forall\,t>0.
\end{equation}
}
As discussed in \cite{Paouris2018}, for any convex function $f$ with $E|f(Z)|^2<\infty,$ the inequality (\ref{a2d3}) implies the following estimates
\begin{equation}\label{faa2}
P\left(f(Z)-M  \leq-t\right)\leq \frac{1}{2}\exp\left(-\frac{\pi t^2}{1024{\rm Var}f(Z) }\right)\,\,\,\forall\,t>0
\end{equation}
and
\begin{align}
P\left(f(Z)-Ef(Z)  \leq-t\right)&\leq \frac{1}{2}\exp\left(-\frac{\pi (t-\sqrt{{\rm Var}f(Z)})^2}{1024{\rm Var}f(Z) }\right)\notag\\
&<\exp\left(-\frac{t^2}{1000{\rm Var}f(Z)}\right)\,\,\,\forall\,t>\sqrt{{\rm Var}f(Z)}.\label{aa2}
\end{align}
%for all $t>1$ where $c>0$ is an absolute constant (from the relation (2.4) in \cite{Paouris2018}, $c=1/1000$).
Interestingly, the bound (\ref{faa2}) and (\ref{aa2}) depend on ${\rm Var}f(Z) $ instead of Lipschitz constant $L.$ Hence, those bound improve (\ref{aa1}) and (\ref{aba1}). Here the improvement lies in the fact that ${\rm Var}f(Z)\leq L^2$ (by the Gaussian Poincar\'e inequality) and one can find many examples of $f$ for which ${\rm Var}f(Z)\ll L^2.$ For further results, Valettas proved in \cite{Valettas2019} that (\ref{aa2}) is tight if the convex function $f$ is not superconcentrated, some variance-sensitive concentration inequalities for the log-concave probability measures were obtained in \cite{Paouris2019,Valettas2019}. We also refer to \cite{Adamczak2005,Samson2003} for other concentration results.

In the present paper, our purpose is to obtain a multi-dimensional version of the variance-sensitive concentration inequality (\ref{aa2}). Our main result is the following statement.
\begin{thm}\label{9jkd4} Let $f_i:\mathbb R^n \to \mathbb R,i=1,\cdots,d$ be twice differentiable functions. For each $i,j,l=1,\cdots,d,$ we define the functions
$$T_{ij}(z):=\int_0^1E\left[\sum\limits_{k=1}^n\frac{\partial  f_i}{\partial z_k}(z)\frac{\partial f_j}{\partial z_k}(\alpha z+\sqrt{1-\alpha^2}Z')\right]d\alpha,\,\,\,z\in \mathbb{R}^n,$$
and
$$T_{lij}(z):=\int_0^1E\left[\sum\limits_{k=1}^n\frac{\partial  f_l}{\partial z_k}(z)\frac{\partial T_{ij}}{\partial z_k}(\alpha z+\sqrt{1-\alpha^2}Z')\right]d\alpha,\,\,\,z\in \mathbb{R}^n.$$%where $Z'$ is an independent copy of $Z.$
Assume that

\noindent (i) $E\left[\exp\left(-\sum\limits_{i=1}^d \theta_i f_i(Z)\right)\right]<\infty$ for any $\theta_1,\cdots,\theta_d\geq 0,$ and all derivatives of $f_i's$ have subexponential growth at infinity,

\noindent (ii) $T_{lij}(z)\geq 0$ for every $i,j,l=1,\cdots,d,$ and $z\in \mathbb{R}^n.$

\noindent Then, for any $t=(t_1,\cdots,t_d)\in \mathbb{R}^d_+,$ we have
\begin{equation}\label{8uiq}
P\left(f_1(Z)-Ef_1(Z)\leq -t_1,\cdots,f_d(Z)-Ef_d(Z)\leq -t_d\right)\leq \exp\left(-\frac{t_1^2+\cdots+t_d^2}{2\|\Sigma\|_{\rm op}}\right),
\end{equation}
where $\Sigma$ is covariance matrix of $(f_1(Z),\cdots,f_d(Z))$ and $\|.\|_{\rm op}$ denotes the operator norm.

Assume, further, that the inverse matrix $\Sigma^{-1}$ exists,  then the bound (\ref{8uiq}) can be improved to the following
\begin{equation}\label{8uiqa}%and  its entries are non-negative
P\left(f_1(Z)-Ef_1(Z)\leq -t_1,\cdots,f_d(Z)-Ef_d(Z)\leq -t_d\right)\leq \exp\left(-\frac{1}{2}t\,\Sigma^{-1}\,t^\intercal\right),
\end{equation}
for any $t$ satisfying $t\,\Sigma^{-1}\in \mathbb{R}^d_+,$ where $t^\intercal$ is the transpose of $t.$
\end{thm}

The proof of Theorem \ref{9jkd4} is deferred to Section \ref{gh3}. To see clearer our new contributions let us end up this Section with some remarks and examples.
\begin{rem}Let $M_i$ be a median for $f _i(Z),i=1,\cdots,d.$ Taking into account the fact that $M_i-Ef_i(Z)\leq E|f_i(Z)-M_i|\leq \sqrt{{\rm Var}f_i(Z)}$ we can get concentration inequalities about the median. Indeed, for example, we obtain from (\ref{8uiq}) that
\begin{align*}
P&\left(f_1(Z)-M_1\leq -t_1,\cdots,f_d(Z)-M_d\leq -t_d\right)\\
&\leq P\left(f_1(Z)-Ef_1(Z)\leq -t_1+M_1-Ef_1(Z),\cdots,f_d(Z)-Ef_d(Z)\leq -t_d+M_d-Ef_d(Z)\right)\\
&\leq P\left(f_1(Z)-Ef_1(Z)\leq -t_1+\sqrt{{\rm Var}f_1(Z)},\cdots,f_d(Z)-Ef_d(Z)\leq -t_d+\sqrt{{\rm Var}f_d(Z)}\right)\\
&\leq \exp\left(-\frac{(t_1-\sqrt{{\rm Var}f_1(Z)})^2+\cdots+(t_d-\sqrt{{\rm Var}f_d(Z)})^2}{2\|\Sigma\|_{\rm op}}\right)
\end{align*}
for all $t_i\geq \sqrt{{\rm Var}f_i(Z)},i=1,\cdots,d.$ Furthermore, by a result of Kwapie\'n \cite{Kwapien1993}, if $f_i's$ are convex functions then $M_i-Ef_i(Z)\leq 0$ and hence,
\begin{align*}
P&\left(f_1(Z)-M_1\leq -t_1,\cdots,f_d(Z)-M_d\leq -t_d\right)\\
&\leq P\left(f_1(Z)-Ef_1(Z)\leq -t_1,\cdots,f_d(Z)-Ef_d(Z)\leq -t_d\right)\\
&\leq \exp\left(-\frac{t_1^2+\cdots+t_d^2}{2\|\Sigma\|_{\rm op}}\right)\,\,\,\forall\,(t_1,\cdots,t_d)\in \mathbb{R}^d_+.
\end{align*}
\end{rem}
\begin{rem} In the one-dimensional case ($d=1$), the bound (\ref{8uiq}) reduces to
\begin{equation}\label{aa2v1}
P\left(f_1(Z)-Ef_1(Z)  \leq-t\right)\leq \exp\left(-\frac{t^2}{2{\rm Var}f_1(Z) }\right)\,\,\,\forall\,\,t\geq 0,
\end{equation}
provided that $T_{111}(z)\geq 0$ for all $z\in \mathbb{R}^n.$ In particular, (\ref{aa2v1}) holds true if $\frac{\partial f_1(z)}{\partial z_i}\geq 0$ and $\frac{\partial^2 f_1(z)}{\partial z_i\partial z_j}\geq 0$ for all $i,j=1,\cdots,n.$ We note that the class of functions $f_1$ satisfying $T_{111}(z)\geq 0$ is not the same as the class of convex functions, see the examples provided by Tanguy in \cite[p. 981]{Tanguy2019}. Hence, our bound (\ref{aa2v1}) supplements the bound (\ref{aa2}). %We also emphasize that, to the best of our knowledge, the bound (\ref{aa2v1}) is first to obtain for general non-linear functions.
\end{rem}
\begin{rem}The concentration inequality (\ref{aa2v1}) was proved first in \cite{ntd2019}, which is a draft version of the present paper. We also refer the reader to \cite[Sec. 5]{Tanguy2019} for discussions and a neater proof using semigroup operators.

\end{rem}
\begin{rem} In the  multi-dimensional setting, an important application of Theorem \ref{9jkd4}  is to estimate the tail of the minimal component of non-Gaussian random vectors. Let $h_i:\mathbb R^n \to \mathbb R,i=1,\cdots,d$ be twice differentiable functions. Assume that the covariance matrix $\Sigma$ of $(h_1(Z),\cdots,h_d(Z))$ is invertible and the functions $f_i=-h_i,i=1,\cdots,d$ satisfy the conditions $(i)$-$(ii)$ of Theorem \ref{9jkd4}. Then, for $(\theta_1,\cdots,\theta_d):=\Sigma^{-1}\,(1,\cdots,1)^\intercal,$ the inequality (\ref{8uiqa}) gives us
\begin{align}
P\left(\min\limits_{1\leq i\leq d}(h_i(Z)-Eh_i(Z))>t\right)&=P\left(f_1(Z)-Ef_1(Z)\leq -t,\cdots,f_d(Z)-Ef_d(Z)\leq -t\right)\notag\\
&\leq \exp\left(-\frac{1}{2}t^2(\theta_1+\cdots+\theta_d)\right)\,\,\,\forall\,\,t\geq 0,\label{uujf3}
\end{align}
provided that $\theta_1+\cdots+\theta_d>0.$ Note that (\ref{uujf3}) is similar to the relation (3.17) given in \cite{Chakrabarty2018} where Chakrabarty \& Samorodnitsky establish the precise asymptotic behaviour for the tail of the minimal component of the Gaussian vector.
\end{rem}
\begin{rem} %In general, it is almost impossible to compute $T_{ij}$ and $T_{lij}$ explicitly. Hence, we need to find a simple way to verify the condition $(ii)$ of Theorem \ref{9jkd4}.Hence, a simple way to verify the condition $(ii)$ of Theorem \ref{9jkd4}.
Observe that
$$\frac{\partial T_{ij}}{\partial z_k}(z)=\int_0^1E\left[\sum\limits_{u=1}^n\frac{\partial^2  f_i}{\partial z_k\partial z_u}(z)\frac{\partial f_j}{\partial z_u}(\alpha z+\sqrt{1-\alpha^2}Z')+\sum\limits_{u=1}^n\frac{\partial  f_i}{\partial z_u}(z)\frac{\partial^2 f_j}{\partial z_k\partial z_u}(\alpha z+\sqrt{1-\alpha^2}Z')\right]d\alpha.$$
Hence, we can verify the condition $(ii)$ of Theorem \ref{9jkd4} by checking the following

\noindent$(ii')$ For all $i=1,\cdots,d,$ $f_i$ is coordinatewise non-decreasing (or non-increasing) and
$$\frac{\partial^2  f_i}{\partial z_k\partial z_u}\geq 0\,\,\forall\,\,k,u=1,\cdots,n.$$
For example, the sum of log-normal random variables fulfills the condition $(ii').$
\end{rem}
\begin{rem}Let $U=(U_1,\cdots,U_n)$ be a Gaussian random vector (coordinates are not necessarily independent). It is well known that we can express $U=AZ,$ where $A$ is a matrix such that $AA^\intercal$ is the covariance matrix of $U.$ Hence, Theorem \ref{9jkd4} can be extended to the functions of arbitrary Gaussian random vectors.
\end{rem}
\begin{exam} Consider $h_1(Z):=Z_1+Z_2$ and $h_2(Z):=Z_1-Z_2,$ we have
$$\Sigma=\left(
           \begin{array}{cc}
             2 & 0 \\
             0 & 2 \\
           \end{array}
         \right).
$$
Hence, the bound (\ref{uujf3}) becomes
\begin{equation}\label{hjd4}
P\left(\min\{Z_1+Z_2,Z_1-Z_2\}>t\right)\leq \exp\left(-\frac{t^2}{2}\right)\,\,\,\forall\,\,t\geq 0.
\end{equation}
On the other hand, we observe that $Z_1+Z_2$ and $Z_1-Z_2$ are independent Gaussian random variables. We obtain
\begin{align*}
P\left(\min\{Z_1+Z_2,Z_1-Z_2\}>t\right)&=P(Z_1+Z_2>t)P(Z_1-Z_2>t)\\
&=(1-\Phi(t/\sqrt{2}))^2\,\,\,\forall\,\,t\geq 0,
\end{align*}
where $\Phi$ is the cumulative distribution function of a standard normal random variable. Thus, in this simple example, our estimate (\ref{hjd4}) is sharp.% on the logarithmic scale.
\end{exam}
\begin{exam} In this example, we construct a class of random vectors satisfying Theorem \ref{9jkd4}. Fixed $N\in \mathbb{N}_+=\{1,2,\cdots\},$ define the set $\mathcal{P}$ of the polynomials of the form
$$\mathcal{P}:=\bigg\{p(x)=\sum\limits_{i=1}^N a_ix^{2m_i},x\in \mathbb{R}:a_i\in\mathbb{R}_+\,\,\,\text{and}\,\,\,m_i\in \mathbb{N}_+,i=1,\cdots,N\bigg\}$$
and denote by $\mathcal{Q}$ the following set of non-negative random variables
$$\mathcal{Q}:=\bigg\{\sum\limits_{i=1}^n p_i(Z_i):p_i\in\mathcal{P}\bigg\}.$$
We observe that, for any $m_1,m_2\in  \mathbb{N}_+,$ the function
$$z_k\mapsto \int_0^1E\left[\sum\limits_{k=1}^nz_k^{2m_1-1}(\alpha z_k+\sqrt{1-\alpha^2}Z_k')^{2m_2-1}\right]d\alpha$$
belongs to $\mathcal{P}$ (this is due to the fact that $E[(Z_k')^{2m-1}]=0$ for any $m\in \mathbb{N}_+$). Hence, it is easy to see that if the random variables $f_i(Z),f_j(Z)\in \mathcal{Q},$ then $T_{ij}(Z)$ also belongs to $\mathcal{Q}.$ As a consequence, we can conclude that every vector $(f_1(Z),\cdots,f_d(Z))\in \mathcal{Q}^d$  will satisfy both conditions $(i)$ and $(ii)$ of Theorem \ref{9jkd4}.

%A very special case of this example is when $f_1(Z)=a_1Z_1^2+...+a_nZ_n^2$ with $a_1,...,a_n\in\mathbb{R}_+.$ In view of (\ref{aa2v1}) we obtain
%\begin{equation}\label{9jfr5}
%P\left(a_1(Z_1^2-1)+...+a_n(Z_n^2-1)\leq-t\right)\leq \exp\left(-\frac{t^2}{4(a_1^2+...+a_n^2) }\right),\,\,t\geq 0.
%\end{equation}
%The above bound recovers Lemma 4.1 in \cite{Li2001}. We also refer the reader to Theorem 4.4 in \cite{Li2001} for an application of (\ref{9jfr5}) in studying of small ball probabilities.

%The Beta distribution
%$$f(Z_1,Z_2)=e^{-Z_1^2-Z_2^2}$$
\end{exam}
\section{Proofs}\label{gh3}
The key tool in our proof is the covariance identify formula for Gaussian functionals. We have
\begin{lem}\label{uuif}Let $f,g:\mathbb R^n \to \mathbb R$ be differentiable functions.  Suppose that $f,g$ and their derivatives have subexponential growth at infinity. Define
$$u(z):=\int_0^1E\left[\sum\limits_{k=1}^n\frac{\partial f}{\partial z_k}(z)\frac{\partial g}{\partial z_k}(\alpha z+\sqrt{1-\alpha^2}Z')\right]d\alpha,\,\,\,z\in \mathbb{R}^n.$$%where $Z'$ is an independent copy of $Z.$
Then, we have
$${\rm Cov}(f(Z),g(Z))= E[u(Z)].$$
\end{lem}
\begin{proof}
Note that the subexponential growth condition ensures the existence of all expectations. Hence, the desired result follows directly from Lemma 2.1.4 in \cite{Adler2007} and the routine approximation argument.
\end{proof}
For the proof of Theorem \ref{9jkd4}, we need the following technical lemma.
\begin{lem}\label{iofh} Let $f_i:\mathbb R^n \to \mathbb R,i=1,\cdots,d$ be twice differentiable functions satisfying the conditions $(i)$-$(ii)$ of Theorem \ref{9jkd4}. Then, for any $\theta_1,\cdots,\theta_d\geq 0,$ we have
\begin{equation}\label{es3r}
E\left[\exp\left(-\sum\limits_{i=1}^d \theta_i f_i(Z)\right)\right]\leq E\left[\exp\left(\sum\limits_{i=1}^d \theta_i N_i\right)\right],
\end{equation}
where $N=(N_1,\cdots,N_d)$ is a Gaussian vector satisfying $E[N_i]=E[f_i(Z)]$ and $E[N_iN_j]=E[f_i(Z)f_j(Z)]$ for all $i,j=1,\cdots,d.$
%$N$ is independent of $(f_1(Z),\cdots,f_d(Z))$  and has the same covariance matrix as that of $(f_1(Z),\cdots,f_d(Z)).$
\end{lem}
\begin{proof}Without loss of generality, we can assume that $N$ is independent of $(f_1(Z),\cdots,f_d(Z))$ and $E[N_i]=E[f_i(Z)]=0$ for all $i=1,\cdots,d.$ Consider the function
$$\varphi(t):=E\left[\exp\left(-\sqrt{t}\sum\limits_{i=1}^d \theta_i f_i(Z)+\sqrt{1-t}\sum\limits_{i=1}^d \theta_i N_i\right)\right],\,\,t\in [0,1].$$
This function is well-defined because of the condition $(i)$ of Theorem \ref{9jkd4}. We have, for $t\in(0,1),$
$$\varphi'(t)=E\left[\exp\left(-\sqrt{t}\sum\limits_{i=1}^d \theta_i f_i(Z)+\sqrt{1-t}\sum\limits_{i=1}^d \theta_i N_i\right)\left(-\frac{1}{2\sqrt{t}}\sum\limits_{i=1}^d \theta_i f_i(Z)-\frac{1}{2\sqrt{1-t}}\sum\limits_{i=1}^d \theta_i N_i\right)\right].$$
By the independence, we obtain
\begin{align}
\varphi'(t)&=-\frac{1}{2\sqrt{t}}E\left[\exp\left(\sqrt{1-t}\sum\limits_{i=1}^d \theta_i N_i\right)\right]E\left[\exp\left(-\sqrt{t}\sum\limits_{i=1}^d \theta_i f_i(Z)\right)\sum\limits_{i=1}^d \theta_i f_i(Z)\right]\notag\\
&-\frac{1}{2\sqrt{1-t}}E\left[\exp\left(-\sqrt{t}\sum\limits_{i=1}^d \theta_i f_i(Z)\right)\right]E\left[\exp\left(\sqrt{1-t}\sum\limits_{i=1}^d \theta_i N_i\right)\sum\limits_{i=1}^d \theta_i N_i\right].\label{p3}
\end{align}
An application of Lemma \ref{uuif} yields
\begin{equation}\label{p2}
E\left[\exp\left(-\sqrt{t}\sum\limits_{i=1}^d \theta_i f_i(Z)\right)\sum\limits_{i=1}^d \theta_i f_i(Z)\right]=-\sqrt{t}E\left[\exp\left(-\sqrt{t}\sum\limits_{i=1}^d \theta_i f_i(Z)\right)T(Z)\right],
\end{equation}
where the function $T$ is defined by
\begin{align*}T(z)&:=\int_0^1E\left[\sum\limits_{k=1}^n\left(\sum\limits_{i=1}^d \theta_i\frac{\partial  f_i}{\partial z_k}(z)\right)\left(\sum\limits_{i=1}^d \theta_i\frac{\partial f_i}{\partial z_k}(\alpha z+\sqrt{1-\alpha^2}Z')\right)\right]d\alpha\\
&=\int_0^1E\left[\sum\limits_{i,j=1}^d \theta_i\theta_j\left(\sum\limits_{k=1}^n\frac{\partial  f_i}{\partial z_k}(z)\frac{\partial f_j}{\partial z_k}(\alpha z+\sqrt{1-\alpha^2}Z')\right)\right]d\alpha\\
&=\sum\limits_{i,j=1}^d \theta_i\theta_jT_{ij}(z),\,\,\,z\in \mathbb{R}^n.
\end{align*}
On the other hand, by  the integration by parts formula for Gaussian random variables (see, e.g. Appendix A.6 in \cite{Talagrand2003}) we obtain
\begin{equation}\label{p1}
E\left[\exp\left(\sqrt{1-t}\sum\limits_{i=1}^d \theta_i N_i\right)\sum\limits_{i=1}^d \theta_i N_i\right]=\sqrt{1-t}E\left[\exp\left(\sqrt{1-t}\sum\limits_{i=1}^d \theta_i N_i\right)\right]\sum\limits_{i,j=1}^d \theta_i\theta_j\sigma_{ij},
\end{equation}
where $\sigma_{ij}:={\rm Cov}(N_i,N_j)={\rm Cov}(f_i(Z),f_j(Z)).$ By inserting (\ref{p2}) and (\ref{p1}) into (\ref{p3}) we get
\begin{align*}
\varphi'(t)=\frac{1}{2}E\left[\exp\left(\sqrt{1-t}\sum\limits_{i=1}^d \theta_i N_i\right)\right]E\left[\exp\left(-\sqrt{t}\sum\limits_{i=1}^d \theta_i f_i(Z)\right)\left(T(Z)-\sum\limits_{i,j=1}^d \theta_i\theta_j\sigma_{ij}\right)\right].
\end{align*}
We now observe that ${\rm Cov}(f_i(Z),f_j(Z))=E[T_{ij}(Z)]$ for every $\,\,i,j=1,\cdots,d.$ Hence,
$$E[T(Z)]=E\left[\sum\limits_{i,j=1}^d \theta_i\theta_jT_{ij}(Z)\right]=\sum\limits_{i,j=1}^d \theta_i\theta_j\sigma_{ij}.$$
So, once again, we can apply Lemma \ref{uuif} and we obtain
$$\varphi'(t)=-\frac{\sqrt{t}}{2}E\left[\exp\left(\sqrt{1-t}\sum\limits_{i=1}^d \theta_i N_i\right)\right] E\left[\exp\left(-\sqrt{t}\sum\limits_{i=1}^d \theta_i f_i(Z)\right)\bar{T}(Z)\right],\,\,t\in(0,1),$$
where $\bar{T}$ is defined by
\begin{align*}
\bar{T}(z)&:=\int_0^1E\left[\sum\limits_{k=1}^n\left(\sum\limits_{i=1}^d \theta_i\frac{\partial  f_i}{\partial z_k}(z)\right)\frac{\partial T}{\partial z_k}(\alpha z+\sqrt{1-\alpha^2}Z')\right]d\alpha\\
&=\int_0^1E\left[\sum\limits_{k=1}^n\left(\sum\limits_{i=1}^d \theta_i\frac{\partial  f_i}{\partial z_k}(z)\right)\left(\sum\limits_{i,j=1}^d \theta_i\theta_j\frac{\partial T_{ij}}{\partial z_k}(\alpha z+\sqrt{1-\alpha^2}Z')\right)\right]d\alpha\\
&=\int_0^1E\left[\sum\limits_{k=1}^n\left(\sum\limits_{i,j,l=1}^d \theta_i\theta_j\theta_l\frac{\partial  f_l}{\partial z_k}(z)\frac{\partial T_{ij}}{\partial z_k}(\alpha z+\sqrt{1-\alpha^2}Z')\right)\right]d\alpha\\
&=\sum\limits_{i,j,l=1}^d \theta_i\theta_j\theta_lT_{lij}(z),\,\,\,z\in \mathbb{R}^n.
\end{align*}
Consequently, the condition $(ii)$ of Theorem \ref{9jkd4} implies that $\bar{T}(z)\geq 0,$ and hence, $\varphi'(t)\leq 0,\,\,t\in(0,1).$ So the claim (\ref{es3r}) holds true because
$$
E\left[\exp\left(-\sum\limits_{i=1}^d \theta_i f_i(Z)\right)\right]=\varphi(1)\leq \varphi(0)= E\left[\exp\left(\sum\limits_{i=1}^d \theta_i N_i\right)\right].
$$
This completes the proof of the lemma.
\end{proof}
\begin{rem} The proof of Lemma \ref{iofh} has some similarities with the proof of Slepian's lemma, see e.g. \cite{Chatterjee2005}. The main difference lies in the fact that we applied Lemma \ref{uuif} twice times. This is the key idea allowing us to handle non-Lipschitz functions.
%Slepian's interpolation
\end{rem}
\begin{rem}The left hand side of (\ref{es3r}) is the Laplace transformation of $(f_1(Z),\cdots,f_d(Z)),$ and hence, the estimate  (\ref{es3r}) could be useful for other research problems.
\end{rem}
\noindent{\bf Proof of Theorem \ref{9jkd4}.} Without loss of generality, we may and will assume that $Ef_i(Z)=0$ for every $i=1,\cdots,d.$ Then, by Markov's inequality we have
\begin{align}
P\big(f_1(Z)-Ef_1(Z)\leq -t_1,\cdots,&f_d(Z)-Ef_d(Z)\leq -t_d\big)\notag\\
&\leq P\left(-\sum\limits_{i=1}^d \theta_i f_i(Z)\geq \sum\limits_{i=1}^d \theta_it_i\right)\notag\\
&\leq E\left[\exp\left(-\sum\limits_{i=1}^d \theta_i f_i(Z)\right)\right]\exp\left(-\sum\limits_{i=1}^d \theta_it_i\right)\label{uuf2}
\end{align}
for all $t=(t_1,\cdots,t_d)\in \mathbb{R}^d_+.$
From (\ref{es3r}) and (\ref{uuf2}) we obtain, for any $\theta=(\theta_1,\cdots,\theta_d)\in \mathbb{R}^d_+,$
\begin{align}
P\big(f_1(Z)-Ef_1(Z)\leq -t_1,\cdots,&f_d(Z)-Ef_d(Z)\leq -t_d\big)\notag\\
&\leq E\left[\exp\left(\sum\limits_{i=1}^d \theta_i N_i\right)\right]\exp\left(-\sum\limits_{i=1}^d \theta_it_i\right)\notag\\
&=\exp\left(-t\,\theta^\intercal+\frac{1}{2}\theta\,\Sigma\, \theta^\intercal\right)\label{ood1}\\
&\leq \exp\left(-t\,\theta^\intercal+\frac{1}{2}\|\Sigma\|_{\rm op}\theta\, \theta^\intercal\right).\label{ood2}
\end{align}
Taking the infimum over all $\theta\in \mathbb{R}^d_+,$ the optimal choice of $\theta$ is $t/\|\Sigma\|_{\rm op}.$ So (\ref{8uiq}) follows from (\ref{ood2}). Similarly, if $\Sigma^{-1}$ exists, we obtain (\ref{8uiqa}) from (\ref{ood1}) by choosing $\theta=t\,\Sigma^{-1}.$

The proof of Theorem \ref{9jkd4} is completed.\hfill$\square$

%\noindent {\bf Acknowledgments.}  The author would like to thank the anonymous referee for their valuable comments and for the references to \cite{Adamczak2005,Chatterjee2005,Kwapien1993,Paouris2019,Pisier1986,Samson2003,Tanguy2019,Valettas2019}.
%This research was funded by Viet Nam National Foundation for Science and Technology Development (NAFOSTED) under grant number 101.03-2019.08.
%\begin{exam}
%$$P\left(\min\{N,1-N^2\}>t\right)$$
%\end{exam}

\end{document}